\RequirePackage{fix-cm}

\documentclass{svjour3}                     % onecolumn (standard format)
\smartqed  % flush right qed marks, e.g. at end of proof
\usepackage{graphicx}

\usepackage{hyperref}
\hypersetup{
    unicode      = false,
    pdftoolbar   = true,
    pdfmenubar   = true,
    pdffitwindow = true,
    pdfnewwindow = true,
    colorlinks   = true,
    linkcolor    = blue,
    citecolor    = blue,
    filecolor    = blue,
    urlcolor     = blue,
    breaklinks   = true
}
\usepackage{graphicx}
\usepackage[english]{babel}
\usepackage[latin1,applemac]{inputenc}  % ligne propre  Apple
\usepackage{amsmath}
\usepackage{pict2e}
\usepackage{amssymb}
\usepackage{tikz}
\usetikzlibrary{positioning,decorations.pathreplacing,shapes}
\usepackage{amsfonts}
\usepackage{multirow}
\usepackage[normalem]{ulem}
\usepackage{epsfig}
\usepackage{soul}
\usepackage{colortbl}
\usepackage{epic}
\usepackage{xspace}
\usepackage{pgf}
\usepackage{rotating}
\usepackage{array}
\usepackage{geometry}
\usepackage{times}
\usepackage[ruled,vlined]{algorithm2e}
\usepackage{algpseudocode}
\usepackage{nomencl}
\usepackage{pgfplots}
\pgfplotsset{width=7cm,compat=1.6}
\usetikzlibrary{patterns}
\usetikzlibrary{arrows}
\newcommand{\argmin}{\mathop{\mathrm{argmin}}}
\newcommand{\argmax}{\mathop{\mathrm{argmax}}}

\newcommand{\poll}{{\sc poll}\xspace}
\newcommand{\search}{{\sc search}\xspace}

\newcommand{\mads}{{\sf Mads}\xspace}

\newcommand{\pb}{{\sf PB}\xspace}
\newcommand{\nomad}{{\sf NOMAD}\xspace}

\def\x{\mathbf{x}}

\newtheorem{Def}{Definition}[section]

\begin{document}

\title{Combining Cross Entropy and MADS methods for \\ inequality constrained global optimization%\thanks{Grants or other notes
}


\author{Charles Audet  	\and
        Jean Bigeon     \and 
        Romain Couderc 
}


\institute{C. Audet \at
            {GERAD}
    		and D\'epartement de math\'ematiques et g\'enie industriel,
    		\'Ecole Polytechnique de Montr\'eal,
    		C.P. 6079, Succ. Centre-ville,
    		Montr\'eal, Qu\'ebec, Canada H3C~3A7. \\
              \email{Charles.Audet@gerad.ca}           %  \\
           \and
           J. Bigeon \at
         	Jean Bigeon, CNRS, LS2N, Nantes, France. \newline
            \email{jean.bigeon@ls2n.fr}
            \and
           R. Couderc \at
            {GERAD}
    		and D\'epartement de math\'ematiques et g\'enie industriel,
    		\'Ecole Polytechnique de Montr\'eal,
    		C.P. 6079, Succ. Centre-ville,
    		Montr\'eal, Qu\'ebec, Canada H3C~3A7 and\\
           Univ. Grenoble Alpes, CNRS, Grenoble INP*, G-SCOP, 38000 Grenoble, France \newline
         *Institute of Engineering Univ. Grenoble Alpes\\
          \email{romain.couderc@grenoble-inp.fr} 
}

\date{Received: date / Accepted: date}

\maketitle

\begin{abstract}
This paper proposes a way to combine the Mesh Adaptive Direct Search (\mads) algorithm with the Cross-Entropy (CE) method for nonsmooth constrained optimization. The CE method is used as an exploration step by the \mads algorithm. The result of this combination retains the convergence properties of \mads and allows an efficient exploration in order to move away from local minima. The CE method samples trial points according to a multivariate normal distribution whose mean and standard deviation are calculated from the best points found so far. Numerical experiments show the efficiency of this method compared to other global optimization heuristics. Moreover, applied on complex engineering test problems, this method allows an important improvement to reach the feasible region and to escape local minima. \\

\end{abstract}

\section*{Declarations}
\begin{itemize}
    \item \textbf{Funding:} Audet is supported by Ivado's fundamental research grant PRF-2019-8079623546. Couderc is supported by a French ministerial grant 8542z.
    \item \textbf{Conflict of interests:} The authors declare that they have no conflict of interest.
    \item \textbf{Availability of data and material:} The authors declare that they do not use any data. 
    \item \textbf{Code availability:} The authors declare that they use the open source nomad software (available at https://www.gerad.ca/nomad/) and custom code.
\end{itemize}
\section{Introduction}
%==============================================================

This work studies inequality constrained blackbox optimization problems of the form:

\begin{equation}
    \underset{\textbf{x} \in \Omega \subseteq \mathbb{R}^n}{ \min} f(\textbf{x}), 
    \label{pb1}
\end{equation}
with
$$\Omega = \{\textbf{x} \in \mathcal{X}: c_j(\textbf{x}) \leq 0, j = 1, 2, ..., m\},$$
where $f:\mathbb{R}^n \rightarrow \Bar{\mathbb{R}} = \mathbb{R}\cup \pm \infty$, $c:\mathbb{R}^n \rightarrow \Bar{\mathbb{R}}^m$ and the set $\mathcal{X}$ represents bound constraints of type $\boldsymbol{\ell} \leq \textbf{x} \leq \textbf{u}$ with $\boldsymbol{\ell}, \textbf{u} \in \Bar{\mathbb{R}}^n$. The specificity of this work is due to the form of the objective function $f$ and of the constraints $c_j$. They can be the result of a simulation of complex physical phenomena. These simulations can take an important amount of time or present some discontinuities and therefore classical optimization methods are difficult to apply. Especially, when the gradient of the objective function and/or of the constraints is not explicitly known, hard to  compute or its estimation is time consuming. This field is called derivative-free optimization (DFO). In the worst case, the gradient does not even exist, which is called blackbox optimization (BBO). \\

Specialized BBO and DFO algorithms have been developed in order to solve this kind of problem. There are two main categories:  model based algorithms \cite{CoScVibook} and direct search algorithms \cite{AuHa2017}. This work deals with direct search algorithms which benefit from theoretical convergence results and adding some modifications may improve their performance. In particular, the Mesh Adaptive Direct Search (\mads) algorithm \cite{AuDe2006} ensures convergence to a point satisfying necessary conditions based on the Clarke calculus \cite{Clar83a}. This theoretical guarantee is a solid basis for blackbox optimization. However, blackbox optimization algorithms must take into account two other types of difficulties. First, algorithms must be efficient in terms of simulation evaluations (constraints and objective function). Indeed, the simulation in an engineering context is often time consuming. Second, blackbox simulations may involve multi-extrema functions, so algorithms must be able to escape from local minima. \mads may be trapped in a local minimum.\\

To address the second difficulty, the \mads algorithm may be combined with Variable Neighborhood Search (VNS) \cite{AuBeLe08} and Latin Hypercube Sampling (LHS) \cite{Stei87a} techniques to escape local minima. Other heuristics of global optimization with no convergence guarantees exist including Evolution Strategy with Covariance Matrix Adaptation (CMA-ES)  \cite{Hansen2006}, Genetic Algorithm (GA) \cite{Goldberg1989}, Differential Evolution (DE) \cite{spde} or Particule Swarm Optimization (PSO) \cite{JKennedy_REberhart_1995}. However, these heuristics often require a large number of function evaluations which is incompatible with the first difficulty. Moreover, no specific mechanism has been developed to enable these methods to deal with inequality constraints. In contrast, some methods  have been developed to address the first problem, by reducing the overall number of simulation evaluation such as: the use of ensembles of surrogate \cite{AuKoLedTa2016} or of quadratic models \cite{CoLed2011} and the integration of the Nelder-Mead (NM) algorithm \cite{G-2017-90}. These different methods improve the efficiency of the \mads algorithm but do not address the difficulty of local optima.\\

The objective of the present research is to propose an alternative strategy: the Cross Entropy (CE) \cite{CE-Ru-Kr} method in hopes of escaping local minima. This method is a trade-off between a more global search and a limited number of blackbox evaluations. It was introduced in 1997, first in a context of rare events in discrete optimization \cite{Ru97} and then adapted to continuous optimization \cite{CE-Ru-Kr-SP}. The main benefit of using CE is that it often converges rapidly to a promising region in the space of variables. However, this method does not benefit from theoretical guarantees, and once it has found a promising region, it requires a large number of simulation evaluations to improve the local accuracy. The two aims of this work are: global exploration with limited number of iterations while preserving the theoretical convergence guarantees. In this purpose, CE is used as a step in the \mads algorithm. The present work proposes a way to include a CE global exploration strategy within the \mads algorithm.\\
	
This paper is divided as follows, Section 2 proposes an overview of the \mads and Cross Entropy methods. Section 3 presents an algorithm combining CE and \mads. Finally, Section 4 shows the main numerical results comparing the proposed method with other \mads type algorithm and state-of-the-art heuristics. Section 5 concludes on future work and on the contributions of this paper.

%==============================================================
\section{Description of \mads and Cross Entropy algorithms}
%==============================================================

This Section describes the \mads and  CE algorithms.

\subsection{The \mads constrained optimization algorithm}
\label{section3}
The present work considers the \mads algorithm with the progressive barrier (\pb)~\cite{AuDe09a}
	to handle inequality constraints
	and with dynamic scaling~\cite{AuLedTr2014}
	to handle the varying magnitudes of the variables. \mads is a direct search algorithm (see algorithm \ref{mads} below). It proceeds iteratively where the blackbox functions are evaluated at some trial points. These points are either accepted as new iterates or rejected, depending on the value of the objective function and of the constraints violations. A key principle of the \mads algorithm is that the candidate points may only be chosen on a discretization of the space of variables called the mesh. This discretization is adaptative and its fineness is controlled by the mesh size vector $\boldsymbol{\delta}^k \in \mathbb{R}_+^n$. In its simplest form, the mesh \cite{AuLedTr2014} is defined as follows:
$$
M^k = V^k + \{diag(\boldsymbol{\delta}^k)\textbf{z} :\textbf{z} \in \mathbb{Z}^n \}
$$
where $V^k$ is the cache containing all points visited by the start of iteration $k$. A point $\textbf{x}$ belongs to $V^k$ if and only if both the objective function and the constraints were evaluated by the start of iteration $k$. The first set $V^0$ may be initialized by the user or by a collection of points generated by LHS for instance.

Each iteration includes three steps. The first is an optional step called the \search where various strategies may be used to explore the space of variables. In practice, the \search accelerates the convergence to an optimum and it may attempt to escape from local minima. The only rules to follow are that the trial points must remain on the current mesh $M^k$ and that the search terminates in finite time. It is in this step of the algorithm that CE was integrated. 

The second step is mandatory and called the \poll. It is the  where the space of variable is locally explored by following strict rules guaranteeing convergence. The \poll is confined to a region delimited by the so called poll size vector $\boldsymbol{\Delta}^k \in \mathbb{R}_{+}^{n}$ which is taken such that $\boldsymbol{\Delta}^k \leq \boldsymbol{\delta}^k$. This region is defined as follows:
\begin{equation*}
    F^k = \{\textbf{x} \in M^k: |\textbf{x}_j - \textbf{x}_j^k| \leq \boldsymbol{\Delta}_j^k, \; \forall j \in [1, n] \}
\end{equation*}
with $\textbf{x}^k$ the current incumbent solution . Then, this step only consist to select a positive spanning set $\mathbb{D}_{\Delta^k}^k$ such that 
\begin{equation*}
    P^k = \{\textbf{x}^k + \boldsymbol{\delta}^k \textbf{d} : \textbf{d} \in \mathbb{D}_{\boldsymbol{\Delta}^k}^k \}
\end{equation*}
is a subset of $F^k$ of extent $\mathbb{D}_{\boldsymbol{\Delta}^k}^k$ and to evaluate the objective functions and the different constraints at these points.

Finally, the last step updates the mesh and the poll size vectors at the end of each iteration. The values of both vectors are reduced when an iteration fails to improve the current solution and they are increased or remain at the same value otherwise. 
To handle the constraints, the \pb is used. This method is based on the constraints violation function \cite{FlLe02a} 
\begin{equation*}
   h(\textbf{x}) := \left\{
            \begin{array}{ll}
                  \displaystyle\sum_{j = 1}^m (\max \{ c_j(\textbf{x}), 0\})^2 & \mbox{ if } \textbf{x} \in \mathcal{X} \subset \mathbb{R}^n\\
                  \infty  & \mbox{ otherwise}.
            \end{array}
        \right. 
\end{equation*}
The constraints violation function value $h(\textbf{x})$ is equal to 0 if and only if the point $\textbf{x}$ belongs to $\Omega$ and is strictly positive otherwise. This function allows to rank any pair of trial points by using the following dominance relation \cite{AuHa2017}.
\begin{Def}
The feasible point $\x \in \Omega$ is said to dominate $\mathbf{y} \in \Omega$ when $f(\x) < f(\mathbf{y})$.
The infeasible point $\x \in \mathcal{X} \setminus \Omega$ is said to dominate $\mathbf{y} \in \mathcal{X} \setminus \Omega$ when $f(\x) \leq f(\mathbf{y})$ and  $h(\x) \leq h(\mathbf{y})$ with at least one strict inequality.
\end{Def}
The PB method approaches an optimal solution by locally exploring around two incumbent solutions. The feasible incumbent solution $\x^{feas} \in \Omega$ and the infeasible incumbent solution $\x^{inf} $ which is the undominated infeasible point with a value of $h$ lower than a threshold called $h_{\max}$. In practice, the threshold $h_{\max}^k$ decreases progressively toward zero without ever reaching it.  Exploring around $\x^{inf}$ may be interesting because it is possible that while the threshold is pushing towards zero a feasible candidate point with a low objective function is generated. Thus, the poll step is applied around these two incumbent solutions. An iteration of the \mads algorithm with the progressive barrier may be of three types:
\begin{itemize}
    \item A dominating iteration occurs when a dominating  trial point with respect to $\x^{inf}$ or $\x^{feas}$ is found. In this case, the threshold is updated to $h_{\max}^{k+1} = h_I^k$.
    \item An improving iteration occurs when it is not a dominating iteration but a trial point  improves the threshold $h_{\max}^k$. In this case, threshold is updated to $h_{\max}^{k+1} = \max \{h(\mathbf{v}) : h(\mathbf{v}) < h_I^k, \mathbf{v} \in V^{k+1} \}$.
    \item An unsuccessful iteration occurs when it is neither a dominating nor an improving iteration. In this case, the threshold is updated to $h_{\max}^{k+1} = h_I^k$.
\end{itemize}
where $h_I^k$ is defined by
\begin{equation*}
    h_I^k = \left \{ \begin{array}{cl}
	h(\x) & \mbox{ for any } \x \in I^k = \{ \displaystyle \argmin_{\x \in U^k} \{f(\x) : 0 < h(\x) \leq h_{\max}^k \} \},  \\
    \infty & \mbox{ if } I^k = \emptyset,
	 \end{array} \right.
\end{equation*}
with $U^k$ the set of infeasible undominated points. Algorithm 1 provides a description of \mads with progressive barrier algorithm, the reader may consult \cite{AuHa2017} for more details, and to \cite{AuDe09a} for a complete presentation.
\begin{algorithm}[htb!]
\caption{The Mesh Adaptive Direct Search algorithm (\mads)} 
{\sf 0. Initialization:}\\
\hspace*{6mm}\begin{tabular}[t]{|lll}
A set of starting point: $V^0 \subset \mathbb{R}^n$\\
An initial poll size vector: $ \boldsymbol{\Delta}^0$\\
The iteration counter: $k \leftarrow 0$\\
The mesh size adjustment parameter: $\tau \in \mathbb{Q} \cap (0,1)$\\
The initial threshold: $h_{\max}^0 = \infty$\\
Define $\boldsymbol{\delta_j^0} = \min\{\boldsymbol{\Delta_j}^0, (\boldsymbol{\Delta_j}^0)^2 \}, \; \forall j \in [1,n]$\\
\end{tabular}

{\sf 1. Search step (optional):}\\
\hspace*{6mm}\begin{tabular}[t]{|lll}
Launch the simulation on a finite set $S^k$ of mesh points.\\
If successful, go to {\sf 3}.
\end{tabular}

{\sf 2. Poll step:}\\
\hspace*{6mm}\begin{tabular}[t]{|lll}
Launch the simulation on the set $P^k$ of poll points.\\
\end{tabular}

{\sf 3. Updates:}\\
\hspace*{6mm}\begin{tabular}{|l}
Update the cache $V^{k+1}$.\\
If the iteration is dominating : \\
\hspace*{6mm} update $\textbf{x}^{k+1}$, $h_{\max}^{k+1} = h_I^k$  and $\boldsymbol{\Delta}^{k+1} = \tau^{-1} \boldsymbol{\Delta}^{k}$.\\
Else if the iteration is improving:\\
\hspace*{6mm} update $\textbf{x}^{k+1}$,$h_{\max}^{k+1} = \max \{h(\mathbf{v}) : h(\mathbf{v}) <h_I^k, \mathbf{v} \in V^{k+1} \}$ and $\boldsymbol{\Delta}^{k+1} =  \boldsymbol{\Delta}^{k}$. \\
Otherwise:\\
\hspace*{6mm} Set $\boldsymbol{\Delta}^{k+1} = \tau \boldsymbol{\Delta}^{k}$ and $h_{\max}^{k+1} = h_I^k$.\\
Set $\boldsymbol{\delta_j^{k+1}} = \min\{\boldsymbol{\Delta_j}^{k+1}, (\boldsymbol{\Delta_j}^{k+1})^2 \}, \; \forall j \in [1,n]$.\\
Increase the iteration counter  $k \leftarrow k+1$ and go to {\sf 1}. \\
\end{tabular}
\label{mads}
\end{algorithm}
%-------------------------------------------------------------------------------------------------------------

The fundamental convergence result~\cite{AuDe09a} of the \mads algorithm with progressive barrier states that if the entire sequence of trial points belongs to a bounded set, 
then there exists an accumulation point $\textbf{x}^*$ such that
the generalized directional derivative $f^\circ(\textbf{x}^*; \textbf{d})$ of Clarke~\cite{Clar83a}
is nonnegative in every hypertangent~\cite{Jahn94a} direction $\textbf{d}$ to the domain $\Omega$ at $\textbf{x}^*$
provided that $\textbf{x}^*$ is feasible. A similar result holds for the constraint violation function $h$ over the set $\mathcal{X}$ in situations where the iterates never approach the feasible region.

\subsection{The Cross Entropy method for continuous optimization}

The Cross Entropy method was introduced by Rubinstein in 1997 in the context of a minimization algorithm for estimating probabilities of rare events \cite{Ru97}. Later, it was modified to solve combinatorial optimization problems \cite{CE-Ru-Kr} and then in 2006 to solve continuous problems \cite{CE-Ru-Kr-SP}. The main idea of this method is as follows. First, each optimization problem is transformed into a rare event estimation problem called associated stochastic problem (ASP). For instance, the deterministic problem (1) is transformed as the minimization of the expectation:  

\begin{equation}
    \min_{\textbf{X} \in \Omega, \gamma \in \mathbb{R}} P(f(\textbf{X}) \leq \gamma) = \min_{\textbf{X} \in \Omega, \gamma \in \mathbb{R}} E(I_{\{f(\textbf{X}) \leq \gamma\}})
\end{equation}

where $\textbf{X}$ is a random vector and $I_{\{f(\textbf{X}\leq \gamma\}}$ is the indicator function. Then, this ASP is tackled  efficiently by an adaptive algorithm. This algorithm constructs a sequence of solutions which converging  to the optimal solution of the ASP. The CE method is composed of two iterative steps:
\begin{itemize}
    \item generation a sample of random data according to a density of probability;
    \item density parameters update thanks to the data sampled to create a new sample in the next iteration.
\end{itemize}

It results that this method  often escapes from local minima.\\

\subsubsection{An introductory example}
For clarity, consider the example from \cite{CE-Ru-Kr-SP} of minimizing the function:
\begin{equation}
    f(x) = -e^{-(x-2)^2} - 0.8e^{-(x+2)^2}, x \in \mathbb{R}.
    \label{3}
\end{equation}
The function $f$ has two local minima and a single global minimum at $x = 2$.
\begin{figure}[htb!]
    \hbox{
    \begin{minipage}[c]{.5\linewidth}
        \centering
        \includegraphics[width = 8 cm]{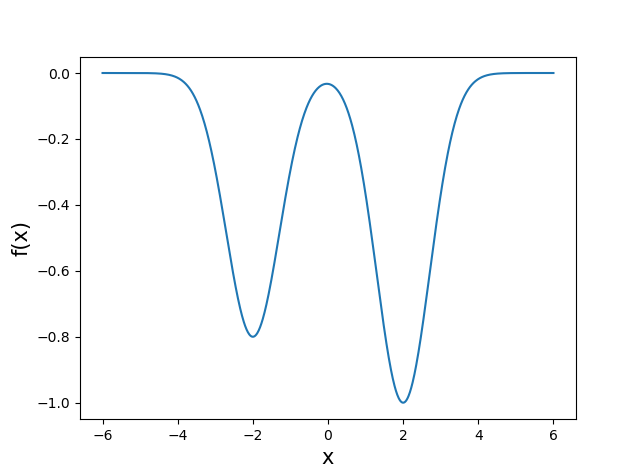}
    \end{minipage}
    \hfill
    \begin{minipage}[c]{.5\linewidth}
        \centering
        \includegraphics[width = 8 cm]{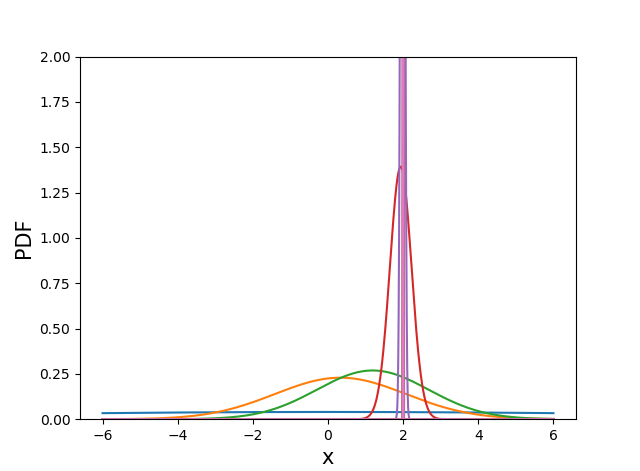}
    \end{minipage}
    }
    \caption{ (Figure inspired by \cite{CE-Ru-Kr-SP}) Graph of objective function $f$ (left) and evolution of the normal distribution during the seven first iterations with $\mu_0 = 0$, $\sigma_0 = 10$,  $N_e =10$ and $N_s = 50$ (right).}
    \label{fig1}
\end{figure}

Using a normal distribution the CE procedure is the following:

\begin{itemize}
    \item Initialization : at the first iteration $k=0$, a mean $\boldsymbol{\mu}^0 \in \mathbb{R}^n$ and a standard deviation $\boldsymbol{\sigma}^0 \in \mathbb{R}^n$ (with $n$ the dimension of the problem) are arbitrarily chosen. A large value of $\boldsymbol{\sigma}^0$ is taken in order to escape from local solutions.
    
    \item Iterative part: at each iteration $k \geq 1$:
    \begin{itemize}
        \item First, a sample $\textbf{X}_1, ..., \textbf{X}_{N_s}$ of points in $\mathbb{R}^n$ is generated from a normal law \\
        $\mathcal{V}(\boldsymbol{\mu}^{k-1},\boldsymbol{\sigma}^{k-1})$ where $N_s$ is the number of samples.
        \item Then, $f$ is evaluated at each sampled points and a number of \textit{elite} points $N_e$, with the lowest value of f. $\boldsymbol{\mu}^k$ and $\boldsymbol{\sigma}^k$ are the mean and standard deviation of these $N_e$ points.
        \item Termination: once the standard deviation becomes sufficiently small, the procedure is stopped. 
    \end{itemize}

\end{itemize}{}
The sequence of normal distribution is illustrated in the right part of Figure \ref{fig1}. This example shows how the CE procedure escapes from the local minimum at $x=-2$ and converges in seven iterations to the  neighborhood of $x=2$.
\subsubsection{The general CE method }

Before presenting the CE method introduced in \cite{CE-Ru-Kr-SP}, the ASP is considered and the two iterative steps of the algorithm are precised. Problem \ref{pb1} is transformed  into an ASP. Using a family of probability distribution functions (pdf) $\{g(\cdot;\textbf{v}): \textbf{v} \in \mathcal{V} \}$  where $g$ is the law chosen to sample the different points at each iteration. $\mathcal{V}$ is the set of vector  parameters of the pdf $g$ which are calculated at each iteration. In the previous example $g$ is taken as the normal law and the $\textbf{v}^k \in \mathcal{V}$ is composed of the mean and standard deviation $\textbf{v}^k = (\boldsymbol{\mu}^k, \boldsymbol{\sigma}^k)$. Having explained the law and its parameters, the ASP related to problem \ref{pb1} can be defined as follows: 
\begin{equation}
    \min_{\textbf{X} \in \Omega, \gamma \in \mathbb{R}} P_{\textbf{v}}(f(\textbf{X}) \leq \gamma) = \min_{\textbf{X} \in \Omega, \gamma \in \mathbb{R}} E_{\textbf{v}}(I_{\{f(\textbf{X}) \leq \gamma\}})
\end{equation}
where $\textbf{v} \in \mathcal{V}$ is a vector of parameter, $\textbf{X}$ is a random vector with a pdf $g(\cdot;\textbf{v})$ and $\gamma$ is a variable. At this stage, for a given value of $\gamma$, the parameter $\textbf{v}$ may be estimated. Conversely, given a vector of parameters $\textbf{v}$, the value $\gamma$ may be also estimated. The CE method  is based on these two estimations, at each iterations, the algorithm estimates one then the other. In the iterative part of Example \ref{3}, the first item corresponds to the estimation of $\gamma$ and the second one to the estimation of $\textbf{v}$. More precisely, we denote $\gamma^* \in \mathbb{R}$ as the infemum of the objective  function, $\textbf{v}^*$ the parameters  and $g(\cdot; \textbf{v}^*)$ the pdf associated to this infemum. The goal is to generate a sequence $(\gamma^k, \textbf{v}^k)$ converging to $(\gamma^*,\textbf{v}^*)$. To achieve this goal, a sequence of pdf $g(\cdot; \textbf{v}^0), g(\cdot;\textbf{v}^1), ...$ converging to $g(\cdot;\textbf{v}^*)$ is created. To assure the convergence, one must have a ``measure" of the difference between the iterate pdf $g(\cdot; \textbf{v}^k)$ and the objective one $g(\cdot;\textbf{v}^*)$. The Kullback-Leibler (KL) divergence \cite{KuLe} is used:
\begin{equation}
    D(g(\cdot;\textbf{v}^*) || g(\cdot; \textbf{v}^k)) = \int_{-\infty}^{\infty}g(\textbf{x};\textbf{v}^*) \ln \left( \frac{g(\textbf{x};\textbf{v}^*)}{g(\textbf{x}; \textbf{v}^k)} \right) d\textbf{x}.
\end{equation}

The iterative steps may now be described. $\rho$ is defined as a very small quantity, corresponding to the proportion of elite points which are kept from an iteration to another. The procedure is:
\begin{itemize}
    \item \textbf{Adaptive update of $\gamma^k$}. With a fixed parameter of pdf $\textbf{v}^{k-1}$, $\gamma^k$ is defined such that it is the $(1-\rho)$-quantile of $f(\textbf{X})$ under $v_{k-1}$. Then, $\gamma^k$ satisfies:
    \begin{equation}
        P_{\textbf{v}^{k-1}}(f(\textbf{X}) \leq \gamma^k) \geq \rho,
        \end{equation}
        \begin{equation}
        P_{\textbf{v}^{k-1}}(f(\textbf{X}) \geq \gamma^k) \geq 1-\rho
        \end{equation}
    where $ X \sim g(\cdot; \textbf{v}^{k-1})$.
    The $\gamma^k$ is denoted $\hat{\gamma}^k$. To obtain this estimator, a sample $\textbf{X}_1, ..., \textbf{X}_{N_s}$ is drawn from $g(\cdot;\textbf{v}^{k-1})$ and evaluated. Then, the $(1-\rho)$ quantile is:
    \begin{equation}
        \hat{\gamma}^k = f_{\lceil (1-\rho)N_s \rceil}.
    \end{equation}
    
    \item \textbf{Adaptive update of $\textbf{v}^k$}. With a fixed $\gamma^k$ and knowing $\textbf{v}^{k-1}$,  $\textbf{v}^k$ is a solution of:
    \begin{equation}
    \begin{split}
                \max_{\textbf{v}} D(\textbf{v}) &= \max_{\textbf{v}} E_{\textbf{v}^{k-1}} \Big( I_{\{f(\textbf{X}) \leq \gamma^k \}} \ln(g(\textbf{X}; \textbf{v})) \Big)\\
                &= \min_{\textbf{v}}E_{\textbf{v}^{k-1}} \Bigg( I_{\{f(\textbf{X}) \leq \gamma^k \}} \ln \left( \frac{I_{\{f(\textbf{X}) \leq \gamma^k \}}}{g(\textbf{X}; \textbf{v})}\right) \Bigg)
    \end{split}
    \end{equation}
    which is the minimization of the KL divergence at iteration $k$ (with the convention $0\ln(0) =0$).
    Nevertheless, in practice, the real expectation and the real $\gamma^k$ are not known, estimators must be used and the following equation is solved:
    \begin{equation}
         \Tilde{\textbf{v}}^k \in \argmax_{\textbf{v}} \widehat{D}(\textbf{v}) =  \frac{1}{N_s} \sum_{i = 1} ^{N_s}  I_{\{f(\textbf{X}_i) \leq \widehat{\gamma}^k \}} \ln(g(\textbf{X}_i; \textbf{v}).
    \end{equation}
\end{itemize}

Last but not least, $\hat{\textbf{v}}^k$ is not set to $\Tilde{\textbf{v}}^k$. There are two reasons for that: first the value of $\Tilde{\textbf{v}}^k$ are smoothed. Second, some components of $\Tilde{\textbf{v}}^k$ could be set to 0 or 1 at the first few iterations and the algorithm could converge to non optimal solution. To avoid these problems, the authors of \cite{CE-Ru-Kr-SP} propose to use the following convex combination:

\begin{equation}
    \hat{\textbf{v}}^k = \alpha \Tilde{\textbf{v}}^k + (1- \alpha) \hat{\textbf{v}}^{k-1}
    \label{eq_sd}
\end{equation}
with $0 < \alpha \leq 1$. Theoretically,  any distribution converging in the neighborhood where the global maximum is attained can be used including normal, double exponential or beta distribution. Nevertheless, the beta distribution has for support $[0,1]$ which is not suitable for global exploration, the double exponential may introduce some discontinuities and the updating step is  quite simple with the normal distribution. Thus, in practice the normal distribution is often chosen \cite{lopez2016gace} \cite{miarnaeimi2018multi}  \cite{6718058}. Therefore, the detailed algorithm is presented next:
\begin{algorithm}[htb!]
\caption{The Cross Entropy (CE) algorithm with normal law} 

Choose $\hat{\boldsymbol{\mu}}^0 \in  \mathbb{R}^n$ and $\hat{\boldsymbol{\sigma}}^0 \in \mathbb{R}^n $\\
Set the iteration counter: $k \leftarrow 0$ \\
$N_s$ number of sampled data at each iteration\\
$N_{e}$ number of elite population\\
$\alpha$ the parameter of convex combination

{\sf 1. Estimation of $\gamma^k$:}\\
\hspace*{6mm}\begin{tabular}[t]{|lll}
Generate a random sample $\textbf{X}_1, .., \textbf{X}_{N_s}$ from $N(\hat{\boldsymbol{\mu}}^{k-1}, \hat{\boldsymbol{\sigma}}^{k-1})$ distribution.\\
Evaluation of the $N_s$ points by the simulation and
then go to {\sf 2}.
\end{tabular}

{\sf 2. Estimation of mean and standard deviation}\\
\hspace*{6mm}\begin{tabular}[t]{|lll}
Let $E^k$ be the indices of the $N_{e}$ best perfoming samples.\\
Set $ \Tilde{\boldsymbol{\mu}}^k = \frac{1}{N_{e}}\displaystyle\sum_{j\in E^k} \textbf{X}_j$ \\
and $[\Tilde{\boldsymbol{\sigma}}^k]_i = \sqrt{\frac{1}{N_{e}-1}\displaystyle\sum_{j\in E^k} [\textbf{X}_j - \boldsymbol{\mu}^k]_i^2}, \; \forall i \in [1, n ]$.\\
\end{tabular}

{\sf 3. Updates:}\\
\hspace*{6mm}\begin{tabular}{|l}
                Apply the convex combinations:\\
                $ \hat{\boldsymbol{\mu}}^k = \alpha \Tilde{\boldsymbol{\mu}}^k + (1- \alpha) \hat{\boldsymbol{\mu}}^{k-1}$\\
                $\hat{\boldsymbol{\sigma}}^k = \alpha \Tilde{\boldsymbol{\sigma}}^k + (1- \alpha) \hat{\boldsymbol{\sigma}}^{k-1}$\\
               Increase the iteration counter  $k \leftarrow k+1$ and go to {\sf 1}. \\
              \end{tabular}

\end{algorithm}

%-------------------------------------------------------------------------------------------------------------

%==============================================================
\section{The CE-MADS constrained optimization algorithm}
\label{sec:MadsNM}
%==============================================================
This section presents the CE-inspired \search step of \mads. Section 4.1 describes how to handle constraints, the update of the mean $\boldsymbol{\mu}$ and the standard deviation $\boldsymbol{\sigma}$ and the condition to enter the CE search step. Section 4.2 presents the algorithm of the CE \search step.

\subsection{The CE-\search step }

\subsubsection{The choice of the elite points}
Section 2 presented the CE method for unconstrained optimization. In \cite{CE-Ru-Kr-SP}, the bound constrained case is treated using a truncated normal law and a penalty approach is used for inequality constraints. In our work, the truncated normal law is also used to treat the bound constraints. For general inequality constraints, the algorithm does not use the penalty approach. The proposed approach is derived from the progressive barrier method described in Section \ref{section3}. In fact, when the algorithm chooses the elite sample, it uses the following function Best (defined in \cite{G-2017-90} and recalled here). Thus, any points in the cache may be selected even if its value of constraint violation is over the threshold of progressive barrier \cite{AuDe09a}. The definition relies on both the objective and the constraint violation functions $f$ and $h$. 

\begin{Def}
\label{def:Best}
The function  $ {\sf Best} : \mathbb{R}^n \times \mathbb{R}^n \mapsto \mathbb{R}^n$
$$ {\sf Best}(\textbf{x},\textbf{y}) = \left \{ \begin{array}{cl}
	\textbf{x} & \mbox{ if } \textbf{x} \mbox{ dominates } \textbf{y} \mbox{ or if } h(\textbf{x}) < h(\textbf{y}),  \\
	\textbf{y} & \mbox{ if } \textbf{y} \mbox{ dominates } \textbf{x} \mbox{ or if } h(\textbf{y}) < h(\textbf{x}),\\
	{\sf Older}(\textbf{x}, \textbf{y}) & \mbox{ Otherwise}
	  \end{array} \right.$$
returns the {\sf best} of two points.
\end{Def}
The function {\sf Older} gives the point which was generated before the former one. Thanks to this definition, CE may treat the general inequality constraints with the terminology used in \mads.

\subsubsection{Update of the mean and standard deviation}

Three elements differ compared to classical CE method concerning the mean and the standard deviation. First, the elite points taken to compute the mean and the standard deviation are not only the $N_s$ points generated by the  normal law. The elite points are chosen from the cache at the iteration $k$, denotes $V^k \subset \mathbb{R}^n$, so any points generated by the \mads algorithm may be selected. This set is ordered with the function {\sf Best}, in order to select the $N_e$ elite points, it is sufficient to take the $N_e$ first points of $V^k$.\\

Second, the mean and the standard deviation initialization procedure differs from the CE method proceeds. Indeed, \mads always begins with a starting point, thus there is at least one point in the cache (the set of evaluated points). Moreover, to avoid generating trial points far from the current solution, bounds are added to the problem as follows (using $\Bar{\textbf{x}}^k$ the poll center at iteration $k$):

\begin{equation}
     (\ell_{i}^k, u_i^k) = \left \{ \begin{array}{cl}
	(\ell_i, u_i) & \mbox{ if } \ell_i \ne -\infty  \mbox{ and } u_i \ne \infty ,\\
	(\Bar{x}_{i}^k   - \tau^{-1}  \Delta_{i}^{\max}, u_i) & \mbox{ if } \ell_i = -\infty  \mbox{ and }  u_i \ne \infty, \\
	(\ell_i, \Bar{x}_{i}^k + \tau^{-1}  \Delta_{i}^{\max}) & \mbox{ if } \ell_i \ne -\infty \mbox{ and } u_i = \infty,\\
	  \end{array} \right. 
	 \;  \forall i \in [1,n].
\end{equation}
where $\Delta^{\max}_i = \max \{\Delta^0_i, ..., \Delta^k_i \}$ for each $i \in [1,n]$. The sequance $(\Delta^{\max}_i)_k$ is non decreasing with respect to $k$ and the product $\tau ^{-1} \Delta_i^{\max}$ is always larger than $\Delta_i$ since $\tau^{-1} > 1$. Once the problem has finite bound constraints, there are two cases to calculate the mean $\boldsymbol{\mu}^k \in \mathbb{R}^n $ and the standard deviation $\boldsymbol{\sigma^k} \in \mathbb{R}^n$:
\begin{itemize}
    \item In case where the number of points in the cache is too small to be relevant, i.e. fewer points that the number $N_e$ required, then the mean and the standard deviation are determined such that:
    \begin{align}
        \boldsymbol{\mu}^k &= \Bar{\textbf{x}}^k \\
        \boldsymbol{\sigma}^k &= 2 (\textbf{u}^k - \boldsymbol{\ell}^k) \label{eq_2}    
    \end{align}
    \item In the others cases, the same calculations are made that in the original CE process: 
    \begin{align*}
        \boldsymbol{\mu}^k &= \frac{1}{N_{e}}\displaystyle\sum_{j\in E^k} \textbf{X}_j \\
        [\boldsymbol{\sigma}^k]_i &=  \displaystyle\sqrt{\frac{1}{N_{e}-1}\displaystyle\sum_{j\in E^k} [\textbf{X}_j - \boldsymbol{\mu}^k]_i^2} , \; \forall i \in [1,n].
    \end{align*}
\end{itemize}

Third, to generate the point during the CE search, the truncated normal law was always used with the bounds created in (12). Moreover, the elite points come not only from the previous normal sampling but also of the other kind of search step. That gives a vector of standard deviation which tends to zero very quickly, the other methods doing generally a local search. That is why,  the standard deviation is calculated as in \ref{eq_sd} with a coefficient $\alpha = 0.7$.

\subsubsection{The condition to pass in the CE-\search step}

The goal of the CE method is to explore in few evaluations the space to determine the promising region. The number of evaluations used by the CE \search step must be quite small. For this purpose, \mads does not perform the CE \search step at every iteration. The standard deviation can be seen as a measure of the incertitude on the data and is used to determine whether to launch the \search step or not. First, a new variable called $\boldsymbol{\sigma}^{p}$ is introduced, it represents the incertitude measured the last time the algorithm passed through the CE \search step and generated trial points. This variable is initialized to $\infty$. Then, the condition to launch the CE \search  is the following:
\begin{equation}
  ||\boldsymbol{\sigma}^{k}||< ||\boldsymbol{\sigma}^{p}|| 
\end{equation}

 This conditions means that the current incertitude is smaller than the previous one. Each time this conditions is respected, $\boldsymbol{\sigma}^p$ is updated with the standard deviation obtained after the CE step.
 Last but not least, there is a special case. The CE method being associate with Mads
which is a local search, it is possible that the points become rapidly close to each others,
reducing the standard deviation. In some cases, that avoids to escape from
unfeasible region. Thus, in case where several iterations of Mads algorithm are passed and the feasible region is still not reached, then the CE- SEARCH is launched with a mean equal to the
current best point and a standard deviation equal to 2 times the initial standard deviation until a feasible point is found.

\subsection{The complete algorithm}

The CE \search step of \mads algorithm is presented here:

\begin{algorithm}[htb!]
\caption{The CE \search step} 
{\sf 1. Calculation of $\boldsymbol{\mu}^k$ and $\boldsymbol{\sigma}^k$ :}\\
\hspace*{6mm}\begin{tabular}[t]{|lll}
if $card(V^k) < N_e$:\\
\hspace*{6mm}\begin{tabular}[t]{|lll}
$\boldsymbol{\mu}^k = \Bar{\textbf{x}}_0$\\
$\boldsymbol{\sigma}^k = 2 (\textbf{u}^k - \boldsymbol{\ell}^k)$ \\
\end{tabular}\\
else:\\
\hspace*{6mm}\begin{tabular}[t]{|lll}
$\boldsymbol{\mu}^k = \frac{1}{N_{e}}\displaystyle\sum_{j\in E^k} \textbf{X}_j$\\
$[\boldsymbol{\sigma}^k]_i =  \sqrt{\frac{1}{N_{e}-1}\displaystyle\sum_{j\in E^k} [\textbf{X}_j - \boldsymbol{\mu}^k]_i^2}, \; \forall i \in [1, n]$ \\
\end{tabular}
\end{tabular}

{\sf 2. CE \search}\\
\hspace*{6mm}\begin{tabular}[t]{|lll}
If $||\sigma^{k}||< ||\sigma^{p}|| $ :\\
\hspace*{6mm}\begin{tabular}[t]{|lll}
Generate a random sample $\textbf{X}_1, .., \textbf{X}_{N_s}$ from $\mathcal{N}(\boldsymbol{\mu}^{k}, 2\boldsymbol{\sigma}^{k})$ distribution \\ and project them on the mesh.\\
 Evaluation of the $N_s$ points by the simulation.\\
 Update: \\
 $\boldsymbol{\mu}^{k+1} = \frac{1}{N_{e}}\displaystyle\sum_{j\in E^k} \textbf{X}_j$\\
$ [\boldsymbol{\sigma}^{k+1}]_i = \sqrt{\frac{1}{N_{e}-1}\displaystyle\sum_{j\in E^k} [\textbf{X}_j - \boldsymbol{\mu}^{k+1}]_i}$ \\
$(\boldsymbol{\sigma}^p)^2 = (\boldsymbol{\sigma}^{k+1})^2$
\end{tabular}
\end{tabular}
\label{algo-Mads}
\end{algorithm}

%================
%==============================================

\section{Computational experiments}

The present work uses data profiles to compare the different algorithm. Data profiles \cite{MoWi2009} allow to assess if algorithms are successful in generating solution values close to the best objective function values. To identify a successful run, a convergence test is required. Let denote $\textbf{x}_e$ the best iterates obtained by one algorithm on one problem after $e$ evaluations, $f_{fea}$ a common reference for a given problem obtained by taking the max feasible objective function values on all run instances of that problem for all algorithms and $f^*$ the best solution obtained by all tested algorithms on all run instances of that problem. Then, the problem is said to be solved within the convergence tolerance $\tau$ when:
$$ f_{fea} - f(\textbf{x}_e) \geq (1-\tau)(f_{fea} - f^*).$$
Different initial points constitute different problems. Moreover, an instance of a problem corresponds to a particular pseudo-random generator seeds.
The horizontal axis of a data profile represents the number of evaluations for problems of fixed dimension, and represents group of $n+1$ evaluations when problems of different dimension are involved. The vertical axis corresponds to the proportion of problems solved within a given tolerance $\tau$. Each algorithm has its curve to allow comparison of algorithms capability to converge to the best objective function value. \\

This section presents the numerical experiments. It is divided in two subsections. The numerical experiments of Section 4.1 are performed on analytical test problems to calibrate the CE-\search parameters. Section 4.2 compares CE-\mads with others state-of-the-art global optimization method. Finally, section 4.3 compares \mads, LH-\mads, VNS-\mads and CE-\mads without the use of models on three real engineering problems.

\subsection{Preliminary experiments to calibrate parameters}

Computational experiments are conducted using the version 3.9.1 of \nomad \cite{Le09b} software package. All tests use the \mads strategy with the use of the NM search \cite{G-2017-90} and without the use of models \cite{AuKoLedTa2016} \cite{CoLed2011}. When the CE-\search is used, it is the first \search step to be applied.

\begin{figure}[htb!]
    \hbox{
    \begin{minipage}[c]{.5\linewidth}
        \centering
        \includegraphics[width = 1 \textwidth]{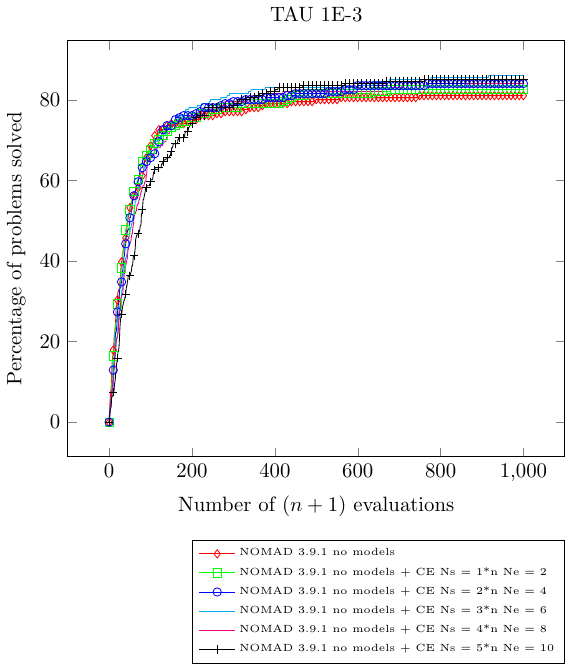}
    \end{minipage}
    \hfill
    \begin{minipage}[c]{.5\linewidth}
        \centering
        \includegraphics[width = 1 \textwidth]{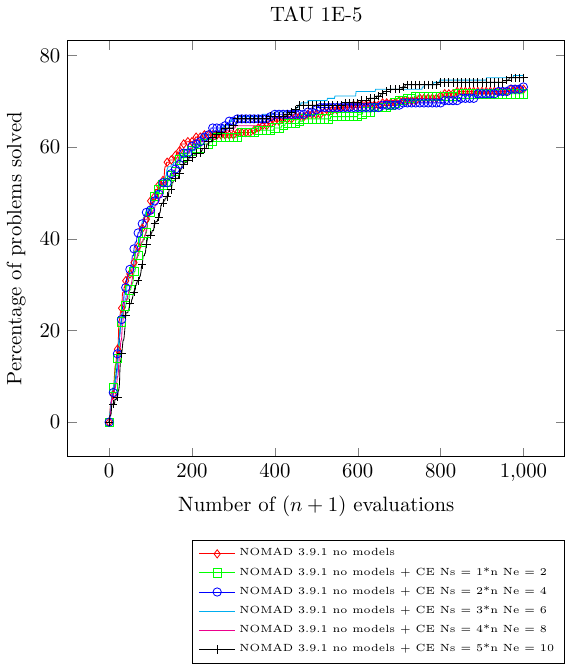}
    \end{minipage}
    }
    \caption{Result of calibration of the hyper-parameters $N_e$ and $N_s$ of CE-MADS on the 69 unconstrained test problems}
    \label{unconstrained}

    \hbox{
    \begin{minipage}[c]{.5\linewidth}
        \centering
        \includegraphics[width = 1 \textwidth]{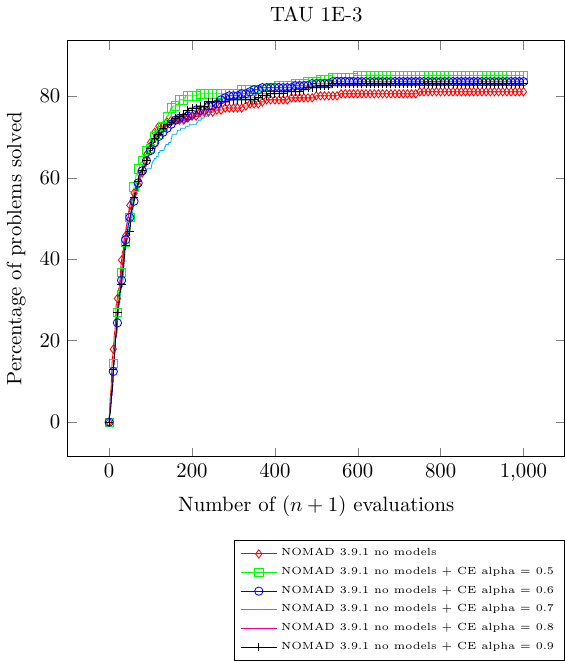}
    \end{minipage}
    \hfill
    \begin{minipage}[c]{.5\linewidth}
        \centering
        \includegraphics[width = 1 \textwidth]{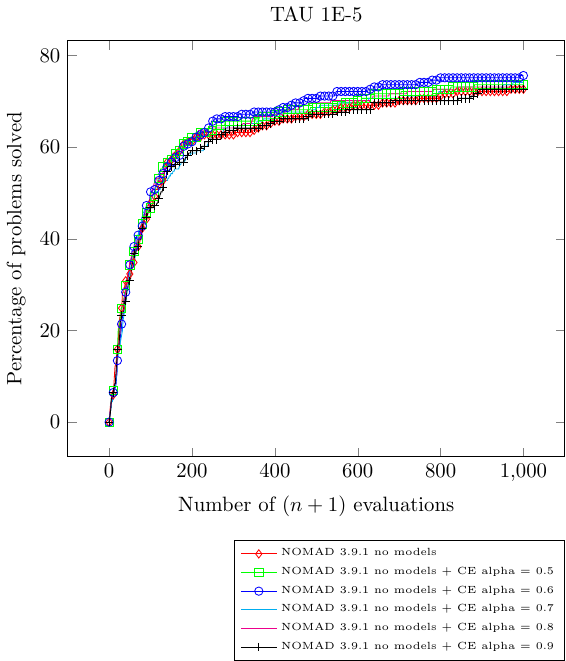}
    \end{minipage}
    }
    \caption{Result of calibration of the hyper-parameter $\alpha$ of CE-MADS on the 69 unconstrained test problems}
    \label{unconstrained_alpha}
\end{figure}

Numerical experiments on analytical test problems are conducted to set default values for the three algorithmic parameters: the parameter of the convex combination $\alpha$, the number of sampled data at each iteration $N_s$ and the number of elite population $N_e$. CE-\mads is tested on 100 analytical problems from the optimization literature. The characteristics and sources of these problems are summarized in Table 1 in appendix \ref{appendix}. The number of variables ranges from 2 to 60; 28 problems have constraints other than bound constraints. In order to have a more precise idea of the effect between the hyper-parameters ($n_e$ and $n_s$), three series of tests are conducted:
\begin{itemize}
    \item A series of tests on the 69 unconstrained test problems having a dimension from 2 to 20.
    \item A series of tests on the 25 constrained test problems having a dimension from 2 to 20.
    \item A series of tests on the 6 larger problems in term of dimension (from 50 to 60), three are constrained and three are not.
\end{itemize}
For each test, the maximal number of function evaluations is set to $1000(n+1)$, where $n$ is the number of variables and each problem is run with $3$ different random seeds.First, for each series of tests, the five following CE-\mads setup of hyper-parameters are compared: $(N_e, N_s) \in \{ (2, n) , (4, 2n), (6, 3n), (8, 4n), (10, 5n) \}$ with $n$ the dimension of the test problem and the $\alpha$ value is fixed to $0.7$ as in the example 3.1 of \cite{CE-Ru-Kr-SP}.A run called \nomad default is added in each series of test to compare our results with the current \nomad software. Data profiles are presented on Figure \ref{unconstrained}, \ref{constrained} and \ref{big} with different values of the tolerance $\tau$.

These results are analysed by series of problems:
\begin{itemize}
    \item On the unconstrained problems (see Figure \ref{unconstrained}), no algorithm really stands out regardless of the value of $\tau$, it is difficult to choose one hyper-parameter rather than another one even if the couple $N_e= 4$ and $N_s = 2n$ appears to be more efficient. 
    \item On the constrained problems (see Figure \ref{constrained}), there are different behaviors according to the value of $\tau$. For $\tau = 10^{-3}$, no algorithm appears to be dominant. However, for $\tau = 10^{-5}$, it happens that greater are the values of $N_s$ and $N_e$, higher is the percentage of problems solved finally. That can be explained because great $N_s$ and $N_e$ allow a better exploration of the space, and so a more precise result at the end.
    \item Finally (see Figure \ref{big}), on the large test problems, and for small values of the tolerance $\tau$ the CE-\mads is outperformed by the \mads algorithm with default values. It seems that the CE method is not useful for problems with a large number of variables. 
\end{itemize}

Second, $N_e$ and $N_s$ are fixed to $4$ and $2n$ repectively Then, the five following CE-\mads setup of hyper-parameters are compared: $\alpha \in \{ 0.5, 0.6, 0.7, 0.8, 0.9 \}$. Data profiles are presented on Figure \ref{unconstrained_alpha}, \ref{constrained_alpha} with different values of the tolerance $\tau$. The tests are not run on the problem with large dimension given the poor performance of CE-\mads on this kind of problems. The results show that none $\alpha$ value outperformed the other ones in any runs of tests. Inspection of the logs of the hyper-parameter calibration reveals the two following observations:
\begin{figure}[htb!]
    \hbox{
    \begin{minipage}[c]{.5\linewidth}
        \centering
        \includegraphics[width = 1 \textwidth]{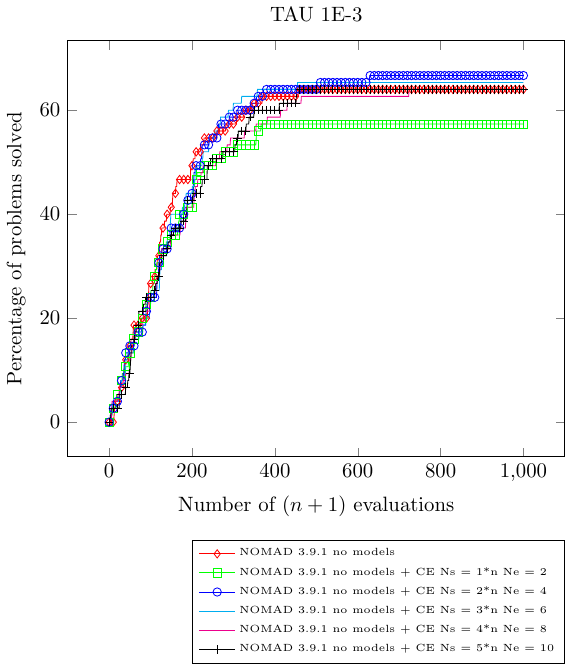}
    \end{minipage}
    \hfill
    \begin{minipage}[c]{.5\linewidth}
        \centering
        \includegraphics[width = 1 \textwidth]{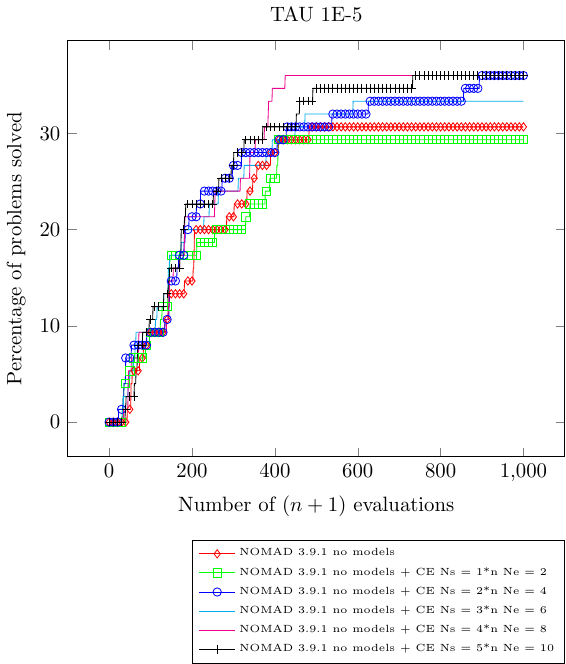}
    \end{minipage}
    }
    \caption{Result of calibration of the hyper-parameters $N_e$ and $N_s$ of CE-MADS on the 25 constrained test problems}
    \label{constrained}

    \hbox{
    \begin{minipage}[c]{.5\linewidth}
        \centering
        \includegraphics[width = 1 \textwidth]{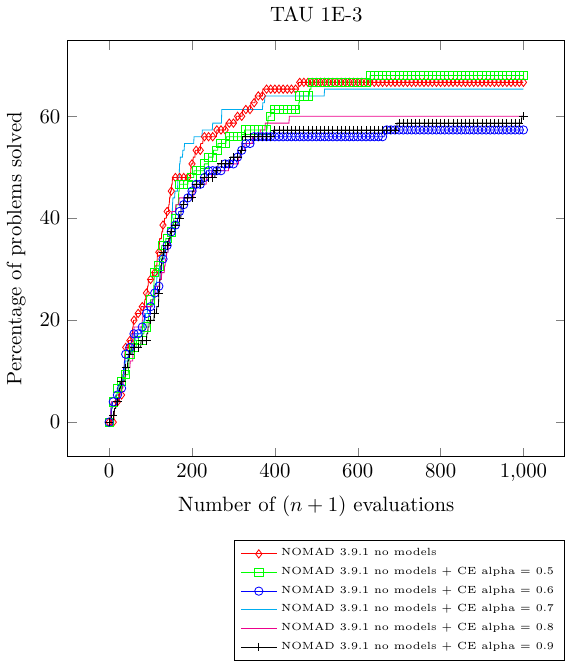}
    \end{minipage}
    \hfill
    \begin{minipage}[c]{.5\linewidth}
        \centering
        \includegraphics[width = 1 \textwidth]{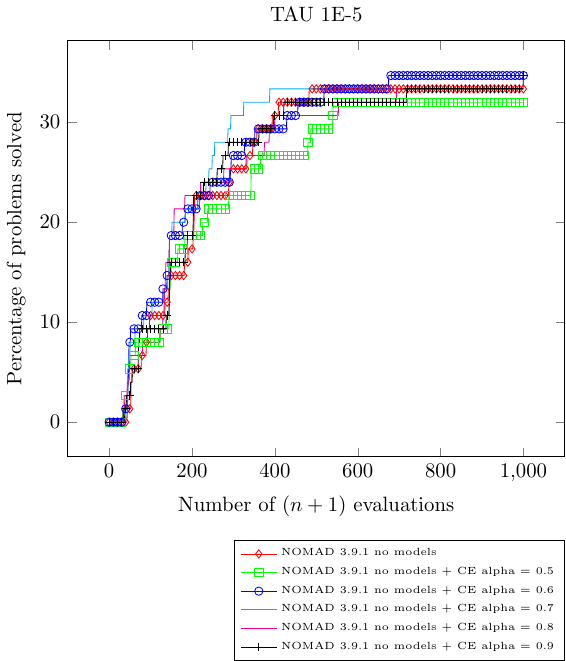}
    \end{minipage}
    }
    \caption{Result of calibration of the hyper-parameter $\alpha$ of CE-MADS on the 25 constrained test problems}
    \label{constrained_alpha}
\end{figure}

\begin{itemize}
    \item The CE-\mads performance is not very sensitive to the hyper-parameter values. This allows to avoid some calibration experiments before applying the algorithm on a new test problem. This is particularly interesting in an engineering context.
    \item For problems with a large number of variables, our tests suggest to avoid using of the CE-\search. Nevertheless, this point has not been confirmed on real engineering problems given that we do not have access to engineering test problems with large dimension.
\end{itemize}

In the remainder of the paper, the CE-\search values are set to $\alpha = 0.7$, $N_e = 4$ and $N_s = 2n$ as they often perform well.

\begin{figure}[htb!]
    \hbox{
    \begin{minipage}[c]{.5\linewidth}
        \centering
        \includegraphics[width = 1 \textwidth]{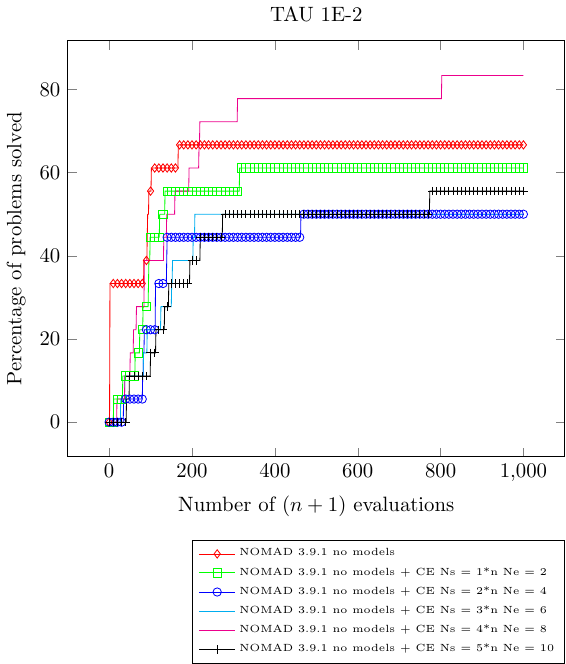}
    \end{minipage}
    \hfill
    \begin{minipage}[c]{.5\linewidth}
        \centering
        \includegraphics[width = 1 \textwidth]{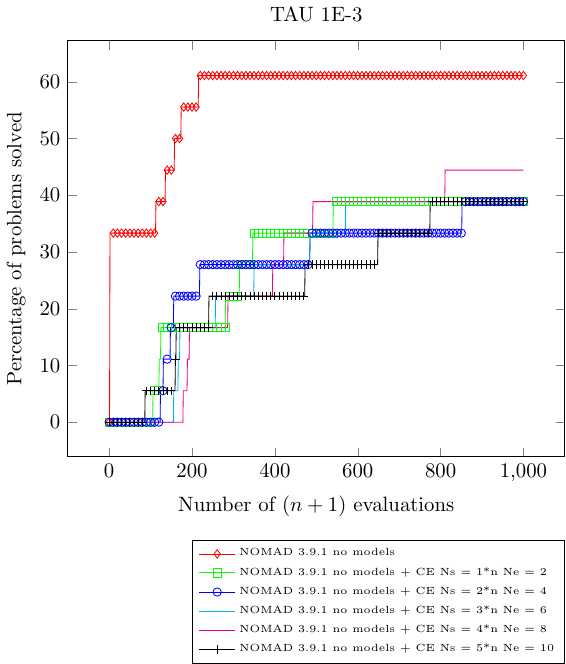}
    \end{minipage}
    }
    \caption{Result of calibration of the hyper-parameters $N_e$ and $N_s$ of CE-MADS on the 6 large test problems}
    \label{big}
\end{figure}

\subsection{Comparisons between CE-\mads and some of state-of-the art global optimization heuristics}

In this section, we compared CE-\mads with other well-known global optimization heuristics designed to escape local minima on a collection of unconstrained global optimization benchmark problems. The framework used to compare the CE-\mads algorithm with others global optimization method is pymoo~\cite{pymoo}. This framework proposes a variety of global optimization algorithm. The CE-\mads method is compared to four of them : 
\begin{itemize}
    \item Genetic Algorithm \cite{Goldberg1989}, a method based on biological inspired operators such as mutation, cross-over and selection. No special advice on the hyper-parameter are given in the pymoo framework, however after few tests, it seems that a population size of 40 performs well. We run the different tests with this value. 
    \item Differential Evolution \cite{price2006differential}, a method which combines evolutionary strategies with geometrical search techniques. In the pymoo framework, the authors advise to test with the following settings: the crossover constant CR  = 0.9, the select weighting factor F = 0.8 and the method is ``DE/rand/1/bin". The size of the population is set to 20 which seems to perform well. 
    
    \item Covariance Matrix Adaptation-Evolutionary Strategy (CMA-ES) \cite{Hansen2006}, a method based also on biological inspired operators. Its name comes from the adaptation of the covariance matrix of the multivariate normal distribution used during the mutation. The setting used for the tests are the default setting.
    
    \item Particle Swarm Optimization \cite{JKennedy_REberhart_1995}, a method inspired by the birds movement and more generally on the collaboration between the individuals. No indication are given in the pymoo framework but it seems that a population size of 15 performs well. 
\end{itemize}

\begin{figure}[htb!]
    \hbox{
    \begin{minipage}[c]{.5\linewidth}
        \centering
        \includegraphics[width = 1.1 \textwidth]{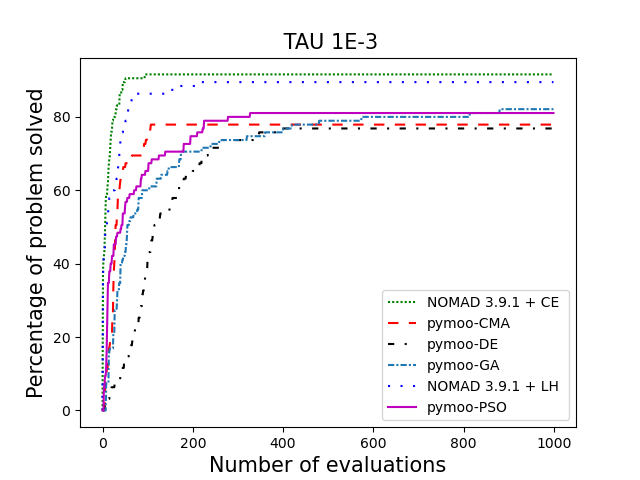}
    \end{minipage}
    \hfill
    \begin{minipage}[c]{.5\linewidth}
        \centering
        \includegraphics[width = 1.1 \textwidth]{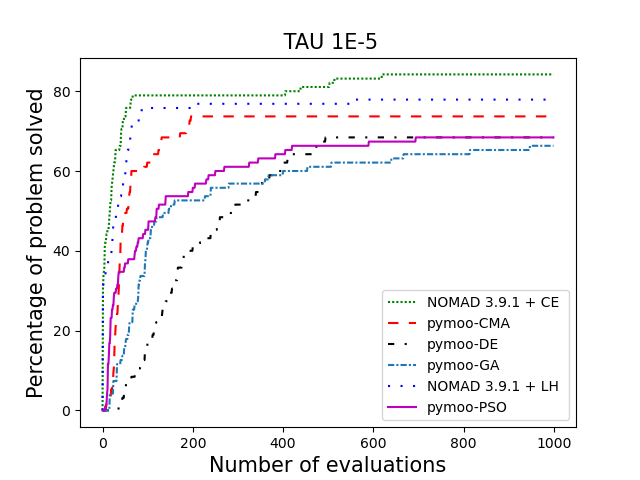}
    \end{minipage}
    }
    \caption{Result on the 19 global optimization test problems  between  CE-\mads, LH-\mads, GA, DE, CMA-ES and PSO for $\tau = 10^{-3}$ (left) and $\tau = 10^{-5}$ (right).}
    \label{global}
\end{figure}

To compare the different algorithms, we use the test problem common between the global optimization benchmark problems of pymoo and the problems provided in the appendix \ref{appendix}. That gives 19 unconstrained test problems (marked in with an asterisk in the appendix) and we run the algorithms with five different seeds in order to reduce the effect of randomness. The maximum number of function evaluations is fixed to 3000. Results are provided on Figure \ref{global}.
The MADS-type algorithms are more efficient than the heuristic of global optimization on this test set. The heuristic's performances are comparable when the required accuracy is $10^{-3}$. If a higher accuracy is desired, it would seem that CMA-ES is the more appropriate method. The gap between heuristics and MADS-type algorithms seems to widen when greater precision is required. This is normal considering that MADS is a local search algorithm originally. The interest of CE-MADS stands out since it allows to combine both a global and a local search, which explains its better performance.

\subsection{Test on engineering problems}
In this section, the CE-\mads algorithm is tested on three different engineering problems and compared  to three algorithms: the \mads-default (without  models), the VNS-\mads where a VNS-\search is used and the LH-\mads which is a default \mads with in addition a LHS search. The comparison with the two last algorithms is crucial because they are methods aiming to explore the space of design variables. The Latin Hypercube \search strategy is used with two parameters $n_{init} = 100$ and $n_{iter} = 10$: $n_{init}$ is the number of LH trial points generated at the first iteration of \mads and $n_{iter}$ the number of LH trial points generated at each subsequent iteration. The Variable Neighbour Search is used with the default parameters \cite{AuBeLe08}. It is an metaheuristic allowing to explore distant neighborhoods of the current incumbent solution. 

\subsubsection{The MDO problem}
\begin{figure}[htb!]
    \hbox{
    \begin{minipage}[c]{.5\linewidth}
        \centering
        \includegraphics[width = 1 \textwidth]{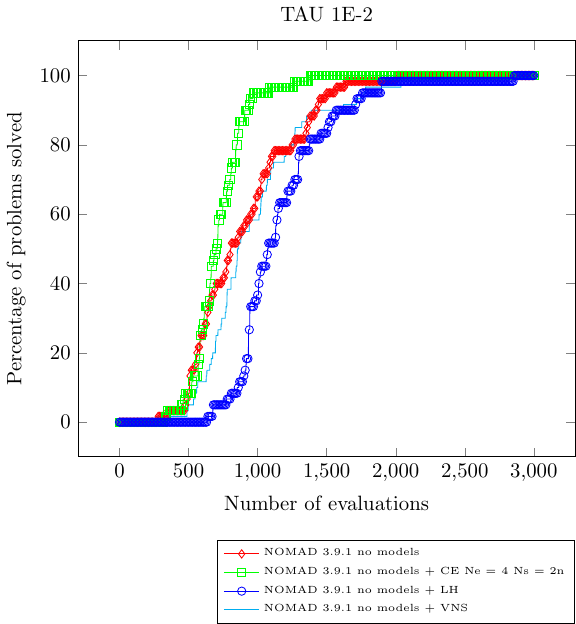}
    \end{minipage}
    \hfill
    \begin{minipage}[c]{.5\linewidth}
        \centering
        \includegraphics[width = 1 \textwidth]{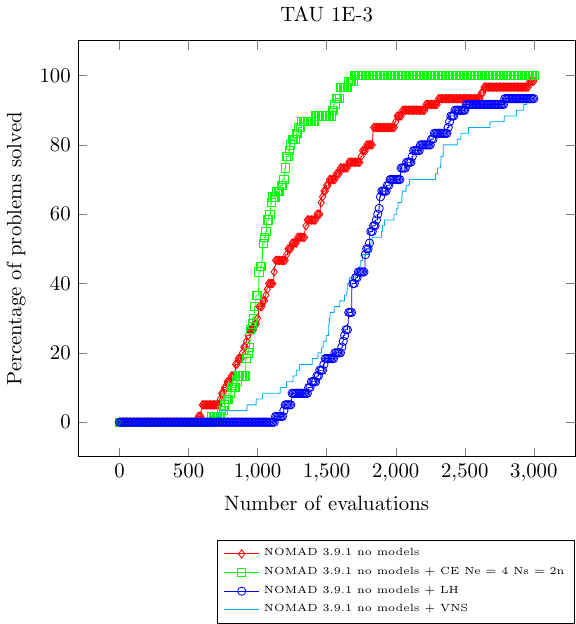}
    \end{minipage}
    }
    \caption{Result on the 60 MDO instances between \mads (no models), CE-\mads, VNS-\mads and LH-\mads for $\tau = 10^{-2}$ (left) and $\tau = 10^{-3}$ (right).}
    \label{mdo}
\end{figure}

The \mads-default (no models), CE-\mads, VNS-\mads and LH-\mads are tested to solve a simple multidisciplinary wing design optimization problem \cite{Giun97a}. Each initial point defines a MDO problem. Solving the problem consists in maximizing the range of an aircraft subject to 10 general constraints. The problem has 10 scaled design variables bounded in $[0; 100]$. Figure \ref{mdo} shows the result on a data profile when solving 20 MDO problems on different initial points using 3000 function evaluations or less. The initial points are real randomly selected within the bounds. Each run is done with three different seeds in order to minimize the impact of the seed. 

\begin{figure}[htb!]
    \hbox{
    \begin{minipage}[c]{.5\linewidth}
        \centering
        \includegraphics[width = 1.1 \textwidth]{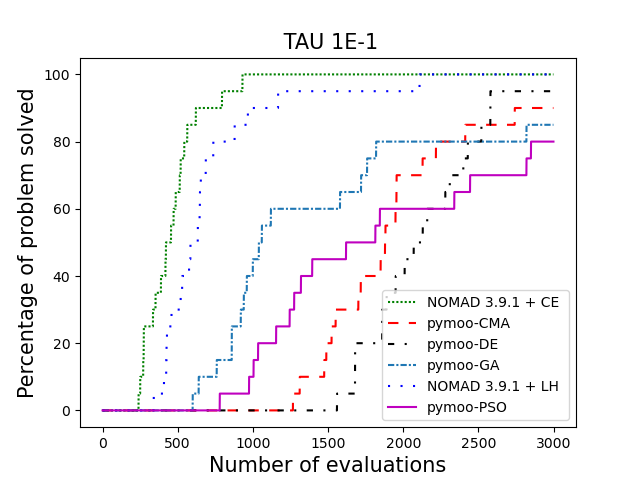}
    \end{minipage}
    \hfill
    \begin{minipage}[c]{.5\linewidth}
        \centering
        \includegraphics[width = 1.1 \textwidth]{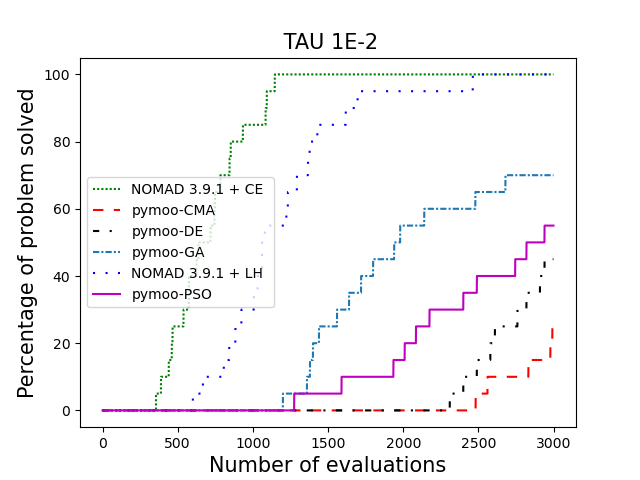}
    \end{minipage}
    }
    \caption{Result on the 20 MDO instances between CE-\mads, LH-\mads, GA, DE, CMA-ES and PSO for $\tau = 10^{-1}$ (left) and $\tau = 10^{-2}$ (right).}
    \label{global2}
\end{figure}

Figure \ref{mdo} shows that the CE-\mads outperforms the other algorithms for all values of $\tau$. Given that the computational time of engineering test problem is relatively low, the comparison between the heuristics and the MADS-type algorithm can be done. The inequality constraints are handled by the default setting in pymoo which is a penalization method. Given that \mads and VNS-\mads require a starting point to be run which can impact their performance, they are not used in this comparison. Therefore, we do not use any starting point but the different algorithms are run with 20 different seeds. The result are given on the figure \ref{global2}. The heuristics perform poorly compared to the \mads-type algorithm. That is not surprising given that \mads benefits of a specialized method to handle the constraints and is specialized in blackbox optimization. The tests on the other engineering test problems are not presented given the great computational time required and the even harder optimization of the blackbox.

\subsubsection{The STYRENE problem}
The \mads-default (no models), CE-\mads, VNS-\mads and LH-\mads algorithms are tested to optimize a styrene production process \cite{AuBeLe08}, called STYRENE. This problem is a simulation of a chemical process. This process relies on a series of interdependent calculation of blocks using common numerical tools as Runge-Kutta, Newton, fixed point and also chemical related solver. The particularity of this problem is the presence of ``hidden" constraints, i.e. sometimes the process does not finish and just return an error. In the case where the chemical process ends, the constraints (not hidden) and the objective functions may be evaluated during a post-processing. The objective is to maximize the net value of the styrene production process with 9 industrial and environmental regulations constraints. 
\begin{figure}[htb!]
    \hbox{
    \begin{minipage}[c]{.5\linewidth}
        \centering
        \includegraphics[width = 1 \textwidth]{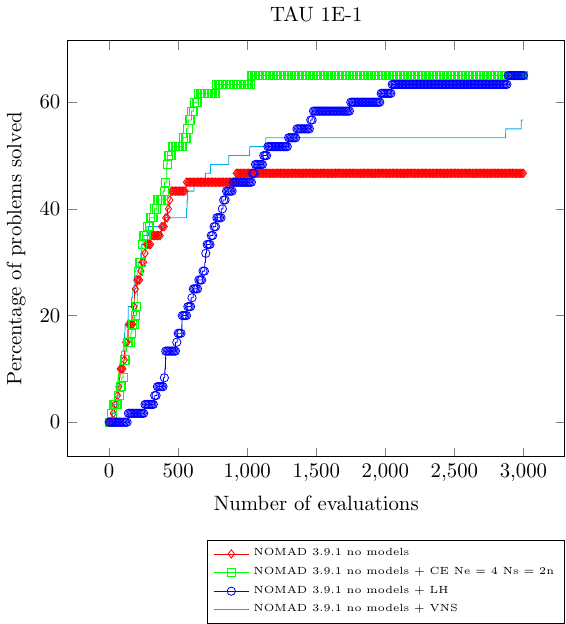}
    \end{minipage}
    \hfill
    \begin{minipage}[c]{.5\linewidth}
        \centering
        \includegraphics[width = 1 \textwidth]{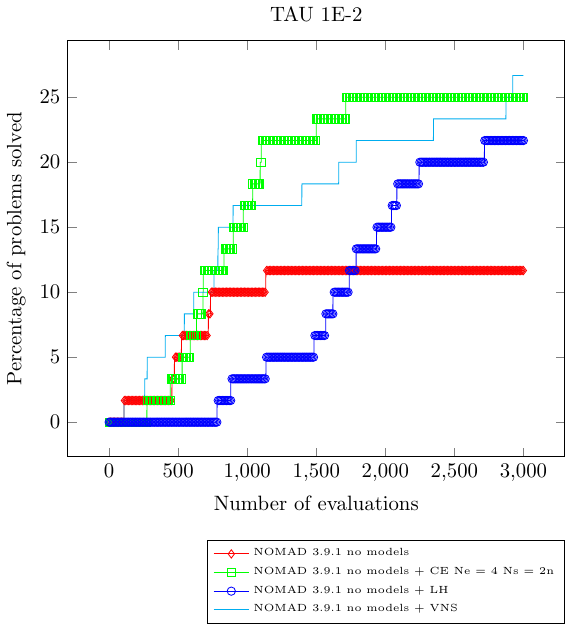}
    \end{minipage}
    }
    \caption{Result on the 60 STYRENE instances between \mads (no models), CE-\mads and LH-\mads for $\tau = 10^{-1}$ (left) and $\tau = 10^{-2}$ (right).}
    \label{Styrene}
\end{figure}
In this work, a STYRENE problem possesses eight independent variables influencing the styrene production process. The variables considered during the optimization process are all scaled and bounded in $X= [0, 100]^8$. As it was done for the MDO test problems, the four algorithms are tested with 20 different starting points taken in $\mathcal{X}$. A maximal number of evaluations of 3000 is used and each problem is run with three different seeds. The STYRENE problems is particularly interesting in this study, because there are two minima as it is shown in \cite{G-2017-90}. The results with $\tau = 10^{-1}$ allow to know the percentage of problems having found the global minimum. The results are provided on Figure \ref{Styrene}. On the left plot, it is interesting to notice that the CE-\mads algorithm find the global minimum the same number of times that the LH-\mads algorithm but is more efficient. On the right plot, the CE-\mads algorithm seems to have the same accuracy that the VNS-\mads algorithm and is slightly more efficient. 

\subsubsection{The LOCKWOOD problem}
\begin{figure}[htb!]
    \hbox{
    \begin{minipage}[c]{.5\linewidth}
        \centering
        \includegraphics[width = 1 \textwidth]{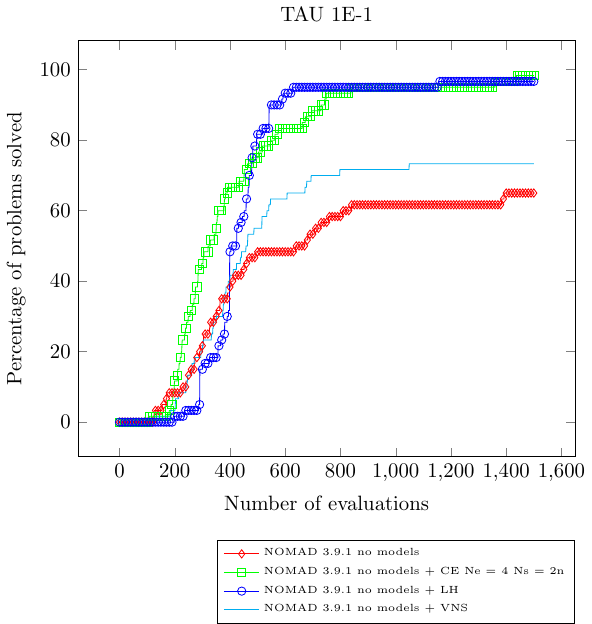}
    \end{minipage}
    \hfill
    \begin{minipage}[c]{.5\linewidth}
        \centering
        \includegraphics[width = 1 \textwidth]{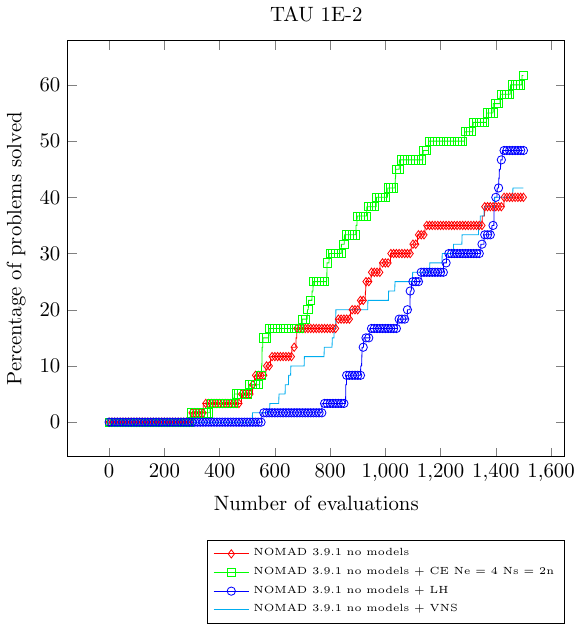}
    \end{minipage}
    }
    \caption{Result on the 20 LOCKWOOD instances between \mads-default (no models), CE-\mads, VNS-\mads and LH-\mads for $\tau = 10^{-1}$ (left) and $\tau = 10^{-2}$ (right).}
    \label{lockwood}
\end{figure}

Finally, the \mads default (no models), LH-\mads and CE-\mads algorithms without quadratic models are tested to solve the basic version of a pump-and-treat groundwater remediation problem from Montana Lockwood Solvent Groundwater Plume Site \cite{matott_lockwood}, called LOCKWOOD. The problem has 6 design variables bounded in $X = [0, 20 000]^6$ and 4 constraints. A particularity of this problem is that each simulation run take several seconds, so the maximum number of blackbox evaluations is set to 1500. The algorithms are started from 20 different randomly selected initial points  in $X$ and three different seeds are used as previously. The results are provided on Figure \ref{lockwood}.
In this problem, reach the feasible region is not easy. Here again, the results at $\tau = 10^{-1}$ allows to give an idea of the number of times the algorithm reach the feasible region. For instance, CE-\mads and LH-\mads always reach the feasible region while \mads default reaches the feasible only 41 times on 60 instances and VNS-\mads only 46 times. The efficiency of CE-\mads and LH-\mads is comparable. However, on the right plot, a better accuracy is reached with a greater efficiency by the CE-\mads algorithm. 

%==============================================================
\section{Discussion}

This paper introduces a way to combine the CE algorithm and the \mads algorithm in order to allow a better space exploration. This is achieved by defining a CE-\search step within the \mads algorithm. The CE-\search generates some points according to a normal distribution whose mean and standard deviation is calculated from the best points stored in the cache. This approach allows to handle the constraints in a different way. Moreover, the particularity of this \search is that it is not performed at each iteration of the \mads algorithm, but according to a criterion based on the value of the norm of the standard deviation of the best points.

Numerical experiments show that in the cases where the problem has different minima or a feasible region hard to reach, the CE-\mads algorithm performs well. Indeed, it attains as often as the LH-\mads the feasible region or the global minimum but it is far more efficient, especially when a tight accuracy is considered. Moreover, even on problem, as MDO, where the classical exploration \search, LH and VNS, do not work well, the CE-\mads algorithm gives interesting results.Finally, comparison with other global algorithms has been made, two conclusions are drawn. First, CE-\mads works better than the heuristics on the unconstrained global optimization test problems. That shows its real ability to escape local optima. Second, CE-\mads outperforms the heuristic on the engineering problems, which is not particularly relevant because it benefits to the \mads ability of performing well on this kind of problems.

Further works will be devoted to improve the link between the \mads algorithm and the CE algorithm by adjusting the size of the mesh with the standard deviation calculated in CE.

\appendix
\clearpage
\section{Appendix}
\label{appendix}

\begin{table}[htb!]
\begin{center}
\renewcommand{\tabcolsep}{3pt}
%\begin{scriptsize}
\begin{footnotesize}
\begin{tabular}{|rrrrrc|}  %cr|}
\hline
\# & Name & Source & $n$ & $m$ & Bnds \\ % & Smth & $f^*$ \\
\hline
\hline
  1  &  ARWHEAD10 	& \cite{GoOrTo03}   & $10$ & $0$ & no  \\[-2 pt]% & yes  & $0.0$ \\
  2  &  ARWHEAD20 	& \cite{GoOrTo03}   & $20$ & $0$ & no  \\[-2 pt] % & yes  & -- \\
  3  &  BARD       		& \cite{MoGaHi81a}                    & $3$ & $0$ & no  \\[-2 pt] % & yes  & -- \\
  4  &  BDQRTIC10 		& \cite{GoOrTo03}     & $10$ & $0$ & no  \\[-2 pt] % & yes  & $18.2812$ \\
  5  &  BDQRTIC20 		& \cite{GoOrTo03}     & $20$ & $0$ & no  \\[-2 pt] % & yes  & -- \\
  6  &  BEALE*         		& \cite{MoGaHi81a}                    & $2$ & $0$ & no  \\[-2 pt] % & yes  & -- \\
  7  &  BIGGS     			& \cite{GoOrTo03}    & $6$  & $0$ & no  \\[-2 pt] % & yes   & $6.97074\cdot10^{-5}$ \\
  8  &  BOX     			& \cite{MoGaHi81a}	  & $3$  & $0$ & no  \\[-2 pt] % & yes   & -- \\ 
  9  &  BRANIN*		& \cite{HeGOTP}       & $2$ & $0$ & yes  \\[-2 pt] % & yes  & $0.397887$ \\
 10 &  BROWNAL5 		& \cite{GoOrTo03}   & $5$ & $0$ & no  \\[-2 pt] % & yes  & -- \\
 11 &  BROWNAL7 		& \cite{GoOrTo03}   & $7$ & $0$ & no  \\[-2 pt] % & yes  & -- \\
 12 &  BROWNAL10 		& \cite{GoOrTo03}   & $10$ & $0$ & no  \\[-2 pt] % & yes  & -- \\
 13 &  BROWNAL20 		& \cite{GoOrTo03}   & $20$ & $0$ & no  \\[-2 pt] % & yes  & -- \\
 14 &  BROWNDENNIS 	    & \cite{MoGaHi81a}   & $4$ & $0$ & no  \\[-2 pt] % & yes  & -- \\ 
 15 &  BROWN\_BS   		& \cite{MoGaHi81a}   & $2$ & $0$ & no  \\[-2 pt] % & yes  & -- \\ 
 16 &  B250   		    & \cite{BoCrFrGaDe07}  & $60$ & $1$ & yes  \\[-2 pt] % & yes  & -- \\ 
 17 &  B500   		    & \cite{BoCrFrGaDe07}   & $60$ & $1$ & yes  \\[-2 pt] % & yes  & -- \\ 
 18 &  CHENWANG\_F2\_X0 & \cite{ChWa2010} & $8$ & $6$ & yes  \\[-2 pt] % &  yes  & -- \\
 19 &  CHENWANG\_F2\_X1	& \cite{ChWa2010} & $8$ & $6$ & yes  \\[-2 pt] % &  yes  & -- \\
 20 &  CHENWANG\_F3\_X0	& \cite{ChWa2010} & $10$ & $8$ & yes \\[-2 pt] %  &  yes  & -- \\
 21 &  CHENWANG\_F3\_X1	& \cite{ChWa2010} & $10$ & $8$ & yes \\[-2 pt] %  &  yes  & -- \\
 22 &  CRESCENT 		& \cite{AuDe09a} & $10$ & $2$ & no   \\[-2 pt] % & yes  & $-9.0$ \\
 23 &  DISK  	     	& \cite{AuDe09a}       & $10$ & $1$ & no  \\[-2 pt] % & yes  & $-17.3205$ \\
 24 &  DIFFICULT2	    & \cite{AuDe09a}       & $10$ & $0$ & no  \\[-2 pt] % & yes  & $-17.3205$ \\
 25 &  ELATTAR 			& \cite{LuVl00}          & $6$ & $0$ & no    \\[-2 pt] % & no  & $0.561139$ \\
 26 &  EVD61 			& \cite{LuVl00}               & $6$ & $0$ & no   \\[-2 pt] %  & no  & $3.51212 \cdot 10^{-2}$ \\
 27 &  FILTER 			& \cite{LuVl00}             & $9$ & $0$ & no    \\[-2 pt] % & no  & $8.40648 \cdot 10^{-2}$ \\
 28 &  FREUDENSTEINROTH* & \cite{MoGaHi81a}  & $2$ & $0$ & no    \\[-2 pt] %  & no  & $8.40648 \cdot 10^{-2}$ \\
 29 &  GAUSSIAN 			& \cite{MoGaHi81a}    & $3$ & $0$ & no    \\[-2 pt] % & no  & -- \\
 30 &  G2\_10				& \cite{AuDeLe07}                & $10$ & $2$ & yes  \\[-2 pt] %  & no  & $-0.740466$ \\
 31 &  G2\_20				& \cite{AuDeLe07}                & $20$ & $2$  & yes \\[-2 pt] % & yes & no  & -- \\
 32 &  G2\_50				& \cite{AuDeLe07}                & $50$ & $2$  & yes \\[-2 pt] % & yes & no  & -- \\
 33 &  GOFFIN			&   \cite{LuVl00}   & $50$ & $0$ & no  \\[-2 pt] %  & --  & -- \\
 34 &  GRIEWANK*			& \cite{HeGOTP}     & $10$ & $0$ & yes  \\[-2 pt] %  & yes  & $0.0$ \\
 35 &  GULFRD*		&   \cite{G-2018-16}   & $3$ & $0$ & no  \\[-2 pt] %  & yes  & $0.0$ \\
 36 &  HELICALVALLEY*	& \cite{MoGaHi81a}     & $3$ & $0$ & no  \\[-2 pt] % & --  & -- \\
 37 &  HS19				& \cite{HoSc1981} & $2$ & $2$ & yes  \\[-2 pt] % &  yes  & -- \\
 38 &  HS78 				& \cite{LuVl00}                  & $5$ & $0$ & no   \\[-2 pt] %   & no  & $-2.49111$ \\
 39 &  HS83\_X0				& \cite{HoSc1981} & $5$ & $6$ & yes  \\[-2 pt] % &  yes  & -- \\
 40 &  HS83\_X1				& \cite{HoSc1981} & $5$ & $6$ & yes  \\[-2 pt] % &  yes  & -- \\
 41 &  HS114\_X0 			& \cite{LuVl00}               & $9$ & $6$ & yes  \\[-2 pt] %  & no  & $-1429.34$ \\
 42 &  HS114\_X1 			& \cite{LuVl00}               & $9$ & $6$ & yes  \\[-2 pt] %  & no  & $-1429.34$ \\
 43 &  JENNRICHSAMPSON & \cite{MoGaHi81a}               & $2$ & $0$ & no  \\[-2 pt] %  & no  & -- \\
 44 &   KOWALIKOSBORNE* & \cite{MoGaHi81a}               & $4$ & $0$ & no  \\[-2 pt] %  & no  & -- \\
 45 &  L1HILB 		&           \cite{LuVl00}    & $50$ & $0$ & no  \\[-2 pt] %  & no  & -- \\
 46 &  MAD6\_X0				& \cite{LuVl00}                & $5$ & $7$ & no  \\[-2 pt] %    & no  & $0.101831$ \\
 47 &  MAD6\_X1			& \cite{LuVl00}                & $5$ & $7$ & no  \\[-2 pt] %    & no  & $0.101831$ \\
 48 &  MCKINNON 		& \cite{McKi98a}                & $2$ & $0$ & no  \\[-2 pt] %    & no  & -- \\
 49 &  MEYER*			& \cite{MoGaHi81a}                & $3$ & $0$ & no   \\[-2 pt] %   & no  & -- \\
 50 &  MEZMONTES		& \cite{MezCoe05} & $2$ & $2$ & yes  \\[-2 pt] % &  yes  & -- \\
 
\hline
\end{tabular}
\begin{tabular}{|rrrrrc|} % cr|}
\hline
\# & Name & Source & $n$ & $m$ & Bnds  \\ % & Smth & $f^*$ \\
\hline
\hline
 51 &  MXHILB 			& \cite{LuVl00}               & $50$ & $0$ & no   \\[-2 pt] %   & no  & -- \\
 52 &  OPTENG\_RBF 		& \cite{KiArYa2011}                & $3$ & $4$ & yes  \\[-2 pt] %  & no  & -- \\
 53 &  OSBORNE1			& \cite{MoGaHi81a} & $5$ & $0$ & no   \\[-2 pt] %  & no  & -- \\
 54 &  OSBORNE2			& \cite{LuVl00} & $11$ & $0$ & no  \\[-2 pt] %  & no  & $9.43876 \cdot 10^{-2}$ \\
 55 &  PBC1 				& \cite{LuVl00} & $5$ & $0$ & no  \\[-2 pt] %  & no  & $8.90604 \cdot 10^{-2}$ \\
 56 &  PENALTY1\_4*		& \cite{GoOrTo03} & $4$ & $0$ & no  \\[-2 pt] % &  yes  & -- \\
 57 &  PENALTY1\_10*	& \cite{GoOrTo03} & $10$ & $0$ & no  \\[-2 pt] % &  yes  & $7.08765\cdot10^{-5}$ \\
 58 &   PENALTY1\_20*		& \cite{GoOrTo03} & $20$ & $0$ & no  \\[-2 pt] % &  yes  & -- \\
 59 &   PENALTY2\_4*		& \cite{GoOrTo03} & $4$ & $0$ & no  \\[-2 pt] % &  yes  & -- \\
 60 &   PENALTY2\_10*		& \cite{GoOrTo03} & $10$ & $0$ & no  \\[-2 pt] % &  yes  & $2.95665 \cdot 10^{-4}$ \\
 61 &   PENALTY2\_20*		& \cite{GoOrTo03} & $20$ & $0$ & no  \\[-2 pt] % &  yes  & -- \\
 62 &  PENTAGON		& \cite{LuVl00} & $6$ & $15$    & no   \\[-2 pt] % &  no  & $-1.85962$ \\
 63 &  PIGACHE\_X00 			& \cite{PigMesNog07} & $4$ & $11$    & yes  \\[-2 pt] %  &  no  & $-1.85962$ \\
 64 &  PIGACHE\_X01			& \cite{PigMesNog07} & $4$ & $11$    & yes  \\[-2 pt] %  &  no  & $-1.85962$ \\
 65 &  POLAK2			& \cite{LuVl00} & $10$ & $0$ & no  \\[-2 pt] %    &  no  & $54.5982$ \\
 66 &  POWELL\_BS 		& \cite{MoGaHi81a} & $2$ & $0$ & no  \\[-2 pt] % &  yes  & --\\
 67 &   POWELLSG4* 		& \cite{GoOrTo03} & $4$ & $0$ & no \\[-2 pt] % &  yes  & $0.0$ \\
 68 &  POWELLSG8 		& \cite{GoOrTo03} & $8$ & $0$ & no \\[-2 pt] % &  yes  & $0.0$ \\
 69 &  POWELLSG12 		& \cite{GoOrTo03} & $12$ & $0$ & no \\[-2 pt] % &  yes  & $0.0$ \\
 70 &  POWELLSG20 		& \cite{GoOrTo03} & $20$ & $0$ & no \\[-2 pt] % &  yes  & $0.0$ \\
 71 &  RADAR 			& \cite{mladenovic2003solving} & $7$ & $0$ & yes \\[-2 pt] % &  yes  & -- \\
 72 &   RANA* 				& \cite{jamil2013literature} & $2$ & $0$ & yes \\[-2 pt] % &  yes  & -- \\
 73 &   RASTRIGIN*		& \cite{HeGOTP} & $2$ & $0$ & yes  \\[-2 pt] %  &  yes  & $0.0$ \\
 74 &  RHEOLOGY 		& \cite{AuHa2017} & $3$ & $0$ & no  \\[-2 pt] %  &  yes  & $0.0$ \\
 75 &   ROSENBROCK*		& \cite{MoGaHi81a} & $2$ & $0$ & yes \\[-2 pt] %  &  yes  & $0.0$ \\
 76 &  SHOR				& \cite{LuVl00} & $5$ & $0$              & no  \\[-2 pt] %   &  no  &   $22.6023$ \\
 77 &  SNAKE 			& \cite{AuDe09a} & $2$ & $2$ & no          \\[-2 pt] %  &  yes  & $0.0$ \\
 78 &  SPRING\_X00 	    & \cite{RodRenWat98} 	& $3$ & $4$ & yes \\[-2 pt] %      &  yes  & -- \\
 79 &  SPRING\_X01 	    & \cite{RodRenWat98} 	& $3$ & $4$ & yes \\[-2 pt] %      &  yes  & -- \\
 80 &  SROSENBR6		& \cite{GoOrTo03} & $ 6$ & $0$ & no \\[-2 pt] % &  yes  & $0.0$ \\
 81 &  SROSENBR8		& \cite{GoOrTo03} & $ 8 $ & $0$ & no \\[-2 pt] % &  yes  & $0.0$ \\
 82 &  SROSENBR10		& \cite{GoOrTo03} & $10$ & $0$ & no \\[-2 pt] % &  yes  & $0.0$ \\
 83 &  SROSENBR20		& \cite{GoOrTo03} & $20$ & $0$ & no \\[-2 pt] % &  yes  & $0.0$ \\
 84 &  TAOWANG\_F2\_X00 & \cite{TaoWan08} & $7$ & $4$ & yes  \\[-2 pt] % &  yes  & -- \\
 85 &  TAOWANG\_F2\_X01	& \cite{TaoWan08} & $7$ & $4$ & yes  \\[-2 pt] % &  yes  &
 86 &   TREFETHEN*	& \cite{jamil2013literature} & $2$ & $0$ & yes \\[-2 pt] % &  yes  & -- \\
 87 &  TRIDIA10 			& \cite{GoOrTo03} & $10$ & $0$ & no \\[-2 pt] %      &  yes  & $0.0$ \\
 88 &  TRIDIA20 			& \cite{GoOrTo03} & $20$ & $0$ & no \\[-2 pt] %      &  yes  & $0.0$ \\
 89 &  TRIGONOMETRIC	& \cite{MoGaHi81a} & $10$ & $0$ & no  \\[-2 pt] %     &  yes  & -- \\
 90 &  VARDIM8			& \cite{GoOrTo03} & $8$ & $0$ & no \\[-2 pt] %    &  yes  & $0.0$ \\
 91 &  VARDIM10			& \cite{GoOrTo03} & $10$ & $0$ & no \\[-2 pt] %    &  yes  & $0.0$ \\
 92 &  VARDIM20			& \cite{GoOrTo03} & $20$ & $0$ & no \\[-2 pt] %    &  yes  & $0.0$ \\
 93 &  WANGWANG\_F3 	& \cite{WanWan10} & $2$ & $0$ & yes  \\[-2 pt] % &  yes  & -- \\
 
 94 &  WATSON9			& \cite{MoGaHi81a} & $9$ & $0$ & no \\[-2 pt] %  &  yes  & $0.0$ \\
 95 &  WATSON12		& \cite{MoGaHi81a} & $12$ & $0$ & yes \\[-2 pt] %  &  yes  & $0.0$ \\
 96 &  WONG1			& \cite{LuVl00} & $7$ & $0$ & no  \\[-2 pt] %         &  no  & $680.707$ \\
 97 &  WONG2			& \cite{LuVl00} & $10$ & $0$ & no \\[-2 pt] %        &  no  & $24.306209$ \\
 98 &  WOODS4			& \cite{GoOrTo03} & $4$ & $0$ & no \\[-2 pt] %  &  yes  & $0.0$ \\
 99 &  WOODS12			& \cite{GoOrTo03} & $12$ & $0$ & no \\[-2 pt] %  &  yes  & $0.0$ \\
 100 &  WOODS20			& \cite{GoOrTo03} & $20$ & $0$ & no  \\[-2 pt] %&  yes  & $0.0$ \\
\hline
\end{tabular}
%\end{scriptsize}
\end{footnotesize}
\end{center}
\caption{Description of the set of 100 analytical problems.}
\label{tab-pbs}
\end{table}

\newpage 
\makenomenclature
\renewcommand{\nompreamble}{The following list describes symbols used within the body of the document. In what follows, if the symbol is bold then it is a vector otherwise it is a scalar.}
 
\nomenclature{$c_j$}{The $j^{th}$ constraint}
\nomenclature{$\Omega$}{The feasible set}
\nomenclature{$\mathcal{X}$}{The bounded constraints set of type $\ell \leq x \leq u$}
\nomenclature{$\boldsymbol{\ell}$}{The lower bound of a decision variable}
\nomenclature{$\mathbf{u}$}{The upper bound of a decision variable}
\nomenclature{$\delta^k$}{The mesh size parameter at iteration $k$}
\nomenclature{$\Delta^k$}{The frame size parameter at iteration $k$}
\nomenclature{$\epsilon$}{The stopping criterion}
\nomenclature{$k$}{The iteration counter}
\nomenclature{$\tau$}{The mesh size adjustment parameter}
\nomenclature{$D$}{A positive spanning set}
\nomenclature{$F^k$}{The frame at iteration $k$}
\nomenclature{$M^k$}{The mesh at iteration $k$}
\nomenclature{$h$}{The measure of constraints violation}
\nomenclature{$I_x$}{Indicator function of x}
\nomenclature{$\mathcal{N}$}{ The normal distribution}
\nomenclature{$\sigma$}{The standard deviation}
\nomenclature{$\mu$}{The mean}
\nomenclature{$\mathcal{V}$}{ Set of parameters of a probability density function}
\nomenclature{$g(\cdot;\cdot)$}{A probability density function}
\nomenclature{$v$}{A parameter of a probability density function}
\nomenclature{$D(u||v)$}{Kullback-Leibler  divergence between u and v}
\nomenclature{$\gamma$}{A parameter to estimate in an associated stochastic problem}
\nomenclature{$\rho$}{A percentage of quantile}
\nomenclature{$N_{e}$}{Number of elite population}
\nomenclature{$N_s$}{Number of sampled data at each iteration}
\nomenclature{$V$}{The cache}
\nomenclature{$n$}{The dimension of a problem}
\nomenclature{$\textbf{X}$}{A random vector}
\nomenclature{$\mathcal{E}$}{The expectation}
\nomenclature{$E^k$}{The set of indices of elite points}
\printnomenclature
%==============================================================

%==============================================================
\clearpage

\bibliographystyle{spmpsci} % basic style, author-year citations
\bibliography{bibliography}


\end{document}